\documentclass[11pt]{article}

\usepackage{amsmath}
\usepackage{amsfonts}
\usepackage{amssymb}
\usepackage{amsthm}

\usepackage[english]{babel}

\usepackage[utf8x]{inputenc} 

\usepackage{graphicx}
\usepackage[autostyle]{csquotes}
\usepackage[linktocpage]{hyperref}
\usepackage{fullpage}
\usepackage{slashed}
\usepackage{mathtools}
\usepackage{indentfirst}
\usepackage{bm}
\usepackage{mathrsfs} 
\usepackage{authblk}
\usepackage{tikz}

\newcommand{\intg}{\mathbb{Z}}
\newcommand{\rational}{\mathbb{Q}}

\newcommand{\complex}{\mathbb{C}}

\newcommand{\lb}{\left(}
\newcommand{\rb}{\right)}

\newcommand{\lac}{\left\{}
\newcommand{\rac}{\right\}}

\newcommand{\RN}[1]{%
  \textup{\uppercase\expandafter{\romannumeral#1}}%
}

\title{Towards a q-series for $osp(2|2n)$ }
\author{John Chae}
\affil{Department of Physics and QMAP, UC Davis, 1 Shields Ave, Davis, CA, 95616, USA \\ {yjchae@ucdavis.edu}}
\date{}  

\begin{document}

\maketitle

\begin{abstract}
A series invariant for a certain class of closed 3-manifolds associated with a type $\RN{1}$ Lie superalgebra $sl(m|n)$ was introduced recently. We find a q-series for the other Lie superalgebra of the same type of the minimum rank.
\end{abstract}

\tableofcontents

\section{Introduction}

Motivated by a prediction for a categorification of the Witten-Reshitikhin-Turaev invariant of a closed oriented 3-manifold~\cite{W1, RT, RT2} in \cite{GPPV, GPV}, a topological invariant q-series for graph 3-manifolds associated with a type $\RN{1}$ Lie superalgebra $sl(m \vert n)$ was introduced recently in \cite{FP}.
This q-series denoted as $\hat{Z}_{ab}$ is labeled by a pair of $Spin^c$ structures of the manifolds as opposed to one label for a q-series ($\hat{Z}_{b}$) corresponding to a classical Lie algebra~\cite{P1}. Another core difference is that the q-series invariant $\hat{Z}_{ab}$ decomposes an extension of the WRT invariant, namely, the CGP invariant based on a Lie superalgebra~\cite{CGP}. In the most general topological setting, the CGP invariant is a topological invariant of a triple consisting of a closed oriented 3-manifold, a link and a certain cohomology class~\footnote{It is denoted $\omega \in H^1 (Y; \intg/2\intg)$, which induces coloring on a surgery link; this link is not part of the triplet~\cite{CGP}.}. The construction of the CGP invariant associated to a Lie (super)algebra involves a new facet: the \textit{modified} quantum dimension, which was first introduced in \cite{GP2} for the quantum groups of Lie superalgebra of type $\RN{1}$ and then for the quantum groups of $sl(2)$ at roots of unity~\cite{GPT}. In general, for a complex Lie superalgebra, the standard quantum dimension vanishes, which makes the WRT invariant of links and 3-manifolds trivial. The modified quantum dimension overcomes the obstacle. In the case of a complex simple Lie algebra, the modified quantum dimension enables to extend the WRT invariant to the above mentioned triplet. Furthermore, for a fixed Lie (super)algebra, the modified dimension can be defined for semisimple and nonsemisimple ribbon categories. Specific examples of the former type were given in \cite{GPT} and \cite{GKP1}. The latter type, which is relevant to this paper, was first dealt with in \cite{GKP2}. And then it was expanded into a Lie superalgebra in \cite{AGP} and \cite{Ha} in which finite dimensional irreducible representations of the (unrolled) quantum groups of $sl(m|n)\, (m\neq n)$ at roots of unity was analyzed (the latter paper focuses on $sl(2|1)$). The CGP invariant contains a variety of information such as the multivariable Alexander polynomial and the ADO polynomial of a link. The relation between the CGP invariant and the WRT invariant for $sl(2;\complex)$ was conjectured in \cite{CGP2} and was proved for certain classes of 3-manifolds. 
\newline

In this paper we analyze $\hat{Z}_{ab}$ associated with $osp(2|2)$ for plumbed 3-manifolds and observe that it is either even or odd power series in the examples. The rest of the paper is organized as follows. In Section 2 we review the ingredients involved in the CGP invariant and in $\hat{Z}_{ab}$ for $sl(m|n)$. In Section 3 we present the formula for homological blocks $\hat{Z}_{ab}$ for $osp(2|2)$. And then in Section 4, we apply the formula to a few examples.
\newline

\noindent \textbf{Acknowledgments.} I am grateful to Bertrand Patureau-Mirand and Pavel Putrov for helpful explanations, as well as to Sergei Gukov for reading a draft of this paper.  

\section{Background}

The twist $\theta$, the S-matrix and the modified quantum dimension $d$ for the quantum group $\mathcal{U}_{h}(g)$ of $g=$ type \RN{1} $= sl(m|n), osp(2|2n)$ are given in \cite{GP2}, where h is a formal variable related to $q=e^{h/2}$:
$$
\theta_V (\lambda) = q^{\left\langle \lambda, \lambda + 2\rho   \right\rangle} 1_V,
$$
where $ \lambda$ is the highest weight of a highest weight $\mathcal{U}_{h}(g)$-module $V$ coloring the a link component and $\rho =\rho_{0} - \rho_{1}$ and $\left\langle - , - \right\rangle $ is the bilinear form (see the appendix A for the details of representation theoretic concepts).
$$
S(V,V^{\prime})= \varphi_{\mu + \rho}(sch(V(\lambda))),
$$
where $\lambda$ and $\mu$ are weights of (irreducible) $\mathcal{U}_{h}(g)$-modules $V$ and $V^{\prime}$, respectively, coloring each link component. Moreover, $sch(V(\lambda))$ is a supercharacter of $V$ and $\varphi_{\beta}$ is a map from a group ring $\intg[\Lambda]$ to $\complex[[h]]$\footnote{ $ e^{\alpha} \mapsto q^{2\left\langle \alpha, \beta \right\rangle}\qquad \text{for any weight}\quad \beta$ }. The modified quantum dimension $d$ is
$$
d(\mu) S(\lambda, \mu) = d(\lambda)  S( \mu, \lambda).
$$
\newline
\indent In the case of $sl(m|n)$, for the unrolled quantum group at roots of unity $\mathcal{U}^{H}_{l}(sl(m|n))\, (m\neq n)$, the above formulas modify~\cite{AGP}:
$$
\theta_V (\lambda; l) = \xi^{\left\langle \lambda, \lambda + \pi  \right\rangle} 1_V\qquad \xi =e^{\frac{i 2\pi}{l}}\quad l \geq m+n-1
$$
\begin{equation*}
\begin{split}
S(V,V^{\prime}; l) & = \varphi_{\lambda + \frac{\pi}{2}}(sch(V(\lambda))),\qquad \pi = 2\rho - 2l\rho_{\bar{0}} \in \Lambda_R\\
                   & = \xi^{2 \left\langle \lambda + \frac{\pi}{2}, \mu + \frac{\pi}{2}  \right\rangle} \prod_{\alpha \in \Delta^{+}_{\bar{0}}} \frac{ \lac l\left\langle \lambda + \frac{\pi}{2}, \alpha  \right\rangle \rac_{\xi}}{ \lac \left\langle \lambda + \frac{\pi}{2}, \alpha  \right\rangle \rac_{\xi}} \prod_{\alpha \in \Delta^{+}_{\bar{1}}} \lac \left\langle \lambda + \frac{\pi}{2}, \alpha  \right\rangle \rac_{\xi}
\end{split}
\end{equation*}
$$
d(\mu; l) S(\lambda, \mu; l) = d(\lambda; l)  S( \mu, \lambda; l).
$$
The second equality in the S-matrix element assumes V is a simple $\mathcal{U}^{H}_{l}(sl(m|n))$-module that is typical having dimension D. The changes of the formula are due to the different pivotal structure of $\mathcal{U}^{H}_{l}(sl(m|n))$. Since the modifications are on the representation theory level and don't seem to involve unique features of $sl(m|n)$, we \textit{assume} that the same modifications occur for $osp(2|2n)$. We next summarize the concepts involved in the homological blocks $\hat{Z}_{ab}$ associated with a Lie superalgebra~\cite{FP}.
\newline

\noindent \underline{Generic graph} In type $\RN{1}$ Lie superalgebra case, plumbing graph used in practice needs to satisfy certain conditions due to the functional form of the edge factor in the integrand of $\hat{Z}$; a plumbing graph is called generic, if
\begin{enumerate}
	\item at least one vertex of a graph has degree $>2$
	
	\item $V|_{deg>2} \neq U \sqcup W$ such that $B^{-1}_{IJ}=0$ for some I $\in V$, J $\in W$,
\end{enumerate}
where V is a set of vertices of a graph.
\newline

\noindent \underline{Good Chamber} For $\hat{Z}$ associated with type $\RN{1}$ Lie superalgebra, there is a notion of a good chamber introduced in \cite{FP}. Existence of a good chamber guarantees that $\hat{Z}$ is a power series with integer coefficients. Specifically, an adjacency matrix $B$ of a generic plumbing graph needs to admit a good chamber. It turns out that $osp(2|2)$ requires positive definite (generic) plumbing graphs, its criteria for good chamber existence shown below are opposite of $sl(2|1)$ criteria in \cite{FP}. Good chamber exists, if there exists a vector $\alpha$ whose components are
$$
\alpha_{I} = \pm 1,\quad I \in V|_{deg\neq 2}
$$
such that
$$
B^{-1}_{IJ}\alpha_{I} \alpha_{J}\geq 0 \quad \forall I \in V|_{deg=1},\quad J \in V|_{deg\neq 2}
$$
$$
B^{-1}_{IJ}\alpha_{I} \alpha_{J} > 0 \quad \forall I,J \in V|_{deg=1},\quad I \neq J 
$$
$$
X_{IJ}\quad \text{is copositive}\qquad X_{IJ}= B^{-1}_{IJ}\alpha_{I} \alpha_{J}\quad I,J \in V|_{deg>2}
$$
are satisfied. Furthermore, a good chamber corresponding to such $\alpha$ is
$$
deg(I) = 1 : \begin{cases}
|y_{I}|^{\alpha_{I}} < 1\\
|z_{I}|^{\alpha_{I}} > 1
\end{cases}
$$
$$
deg(I) > 2 : \bigg| \frac{y_{I}}{z_{I}} \bigg|^{\alpha_{I}} < 1.
$$
\newline

\noindent \underline{Chamber Expansion} The series expansions for the good chambers for $osp(2|2)$ turn out to be the same as that of $sl(2|1)$ in \cite{FP}. We record here the expansions.
\begin{itemize}
	\item $\text{degree(I)}=2+ K >2$
\begin{flushleft}
	$
	\lb \frac{(1-z_{I})(1-y_{I})}{y_{I}-z_{I}} \rb^{K} = \begin{cases}
	(z_{I}-1)^K (1-y_{I}^{-1})^K \sum_{r_{I}=0}^{\infty} \frac{(r_I +1)(r_I +2) \cdots (r_I +K-1)}{(K-1)!} \lb z_{I}/ y_{I} \rb^{r_I}, & |z_{I}| < |y_{I}| \\
	(1 - y_{I})^K (1-z_{I}^{-1})^K \sum_{r_{I}=0}^{\infty} \frac{(r_I +1)(r_I +2) \cdots (r_I +K-1)}{(K-1)!} \lb y_{I}/ z_{I} \rb^{r_I}, & |z_{I}| > |y_{I}| \\
	\end{cases}
	$
\end{flushleft}

\item $\text{degree(I)}=1$
	\begin{flushleft}
	$
	 \frac{y_{I}-z_{I}}{(1-z_{I})(1-y_{I})}  = \begin{cases}
	\sum_{r_{I}=1}^{\infty} y_{I}^{r_I} + \sum_{r_{I}=0}^{\infty} z_{I}^{-r_I}, & |y_{I}| < 1, |z_{I}| > 1 \\
	-\sum_{r_{I}=0}^{\infty} y_{I}^{-r_I} - \sum_{r_{I}=1}^{\infty} z_{I}^{r_I}, & |y_{I}| > 1, |z_{I}| < 1 \\
	\end{cases}
	$
\end{flushleft}
	
\end{itemize}

\section{Result}

For closed oriented plumbed 3-manifolds $Y(\Gamma)$ having $b_1(Y)=0$ and $a,b \in Spin^{c}(Y) \cong H_1(Y)$, $\hat{Z}_{ab}[Y;q]$ associated with $osp(2|2)$ is
\begin{align*}
\hat{Z}^{osp(2|2)}_{ab}[Y;q] & = (-1)^{\pi} \int_{\Omega} \prod_{I\in V} \frac{dz_{I}}{i2\pi z_{I}}\frac{dy_{I}}{i2\pi y_{I}} \lb \frac{y_{I}-z_{I}}{(1-z_{I})(1-y_{I})} \rb^{2-deg(I)} \sum_{\substack{ n\in B Z^L + a \\ m\in B Z^L + b} } q^{-2n B^{-1} m} \prod_{J\in V} z_{J}^{m_{J}} y_{J}^{n_{J}} \\
& \in \rational + q^{\Delta_{ab}}\intg [[q]],\qquad \Delta_{ab} \in \rational\qquad H_1(Y) = \intg^{L} / B \intg^{L}\qquad |q|<1
\end{align*}
where $B=B(\Gamma)$ is an adjacency matrix of a generic plumbing graph $\Gamma$. Moreover, the convergent q-series in a complex unit disc requires $\Gamma$ to be positive definite plumbing graphs inside a complex unit disc. This is opposite of the $sl(m|n)$ case (and for $su(n)$). The edge factor in the integrand turns out to be the same as that of $sl(2|1)$~\footnote{This seems to be no longer true in higher rank cases}. The integration prescription states that one chooses a chamber for an expansion of the integrand using the chamber expansion in the previous section and then picks constant terms in $y_I$ and $z_I$.

\section{Examples}

We apply the above formula to $\intg HS^3$ and $\rational HS^3$. Each $\hat{Z}$ is either even or odd power series and the regularized constants are given by the zeta function $\zeta(s)$ or the Hurwitz zeta function $\zeta(s,x)$ (see \cite{FP} for details).
\newline

\begin{itemize}
	\item $Y=S^3$\quad The adjacency matrix of a generic plumbing graph of $S^3$ (Figure 1 in appendix B) is
	
$$
B(\Gamma)= \begin{pmatrix}
4 & 1 & 1 & 1\\
1 & 1 & 0 & 0\\
1 & 0 & 1 & 0\\
1 & 0 & 0 & 1\\
\end{pmatrix}\qquad
Det B = 1\qquad |H_1(Y)|= |DetB|
$$
The two chambers are 
$$
\alpha = \pm (1,-1,-1,-1).
$$
$$
\hat{Z}^{osp(2|2)}_{00}[Y ;q]= 1+ 2\zeta(-1)+2\zeta(0) -2 q^2-4 q^4-4 q^6-6 q^8-4 q^{10}-8 q^{12}-4 q^{14}-8 q^{16}-6 q^{18} + \cdots  
$$
\newline

\item $Y=\Sigma(2,3,7)$\quad The adjacency matrix of a (generic) plumbing graph (Figure 2 in appendix B) is
$$
B(\Gamma)= \begin{pmatrix}
1 & 1 & 1 & 1\\
1 & 2 & 0 & 0\\
1 & 0 & 3 & 0\\
1 & 0 & 0 & 7\\
\end{pmatrix}\qquad
Det B = 1
$$
The two chambers are 
$$
\alpha = \pm (1,-1,-1,-1).
$$
\begin{align*}
\hat{Z}^{osp(2|2)}_{00}& =   1+ 2\zeta(-1)+2\zeta(0) -2 q^4-2 q^6-4 q^8-2 q^{10}-6 q^{12}-4 q^{14}-6 q^{16}-6 q^{18} -8 q^{20} \\
                              & -4 q^{22}-10 q^{24}-6 q^{26}-8 q^{28} + \cdots 
\end{align*}
\newline

\item $Y=M(-1|\frac{1}{2},\frac{1}{3},\frac{1}{8})$\quad The adjacency matrix of a (generic) plumbing graph (Figure 3 in appendix B) is
$$
B(\Gamma)= \begin{pmatrix}
1 & 1 & 1 & 1\\
1 & 2 & 0 & 0\\
1 & 0 & 3 & 0\\
1 & 0 & 0 & 8\\
\end{pmatrix}\qquad
Det B = 2
$$
The two chambers are 
$$
\alpha = \pm (1,-1,-1,-1).
$$
The independent $\hat{Z}_{ab}$ are
\begin{align*}
\hat{Z}_{00} & =   1 + 2\zeta(-1)+2\zeta(0) -2 q^4-4 q^8-2 q^{10}-4 q^{12}-2 q^{14}-6 q^{16}-2 q^{18}-6 q^{20}-2 q^{22}-8 q^{24}\\
                            & -2 q^{26}-4q^{28}  + \cdots\\
\hat{Z}_{11}& =  -2 q^3-2 q^5-2 q^7-4 q^9-4 q^{11}-2 q^{13}-6 q^{15}-4 q^{17}-4 q^{19}-6 q^{21}-4 q^{23}- 6 q^{27} + \cdots \\
\hat{Z}_{01}& = 2\zeta(-1,1/2) + \zeta(0,1/2) -q^2-q^4-3 q^6-2 q^8-3 q^{10}-4 q^{12}-4 q^{14}-3 q^{16}-5 q^{18}-5 q^{20}\\
                              & -4 q^{22}- 5 q^{24} -4 q^{26}-5q^{28} + \cdots														
\end{align*}
where the labels denote the last components of elements of $H_1(Y)$.
\newline

\item $Y=M(-1|\frac{1}{2},\frac{1}{3},\frac{1}{9})$\quad The adjacency matrix of a (generic) plumbing graph (Figure 4 in appendix B) is
$$
B(\Gamma)= \begin{pmatrix}
1 & 1 & 1 & 1\\
1 & 2 & 0 & 0\\
1 & 0 & 3 & 0\\
1 & 0 & 0 & 9\\
\end{pmatrix}\qquad
Det B = 3
$$
The independent $\hat{Z}^{osp(2|2)}_{ab}[Y]$ are
\begin{align*}
\hat{Z}_{00} & = 1 + 2\zeta(-1)+2\zeta(0) -2 q^6-2 q^8-4 q^{12}-2 q^{14}-2 q^{16}-4 q^{18}-2 q^{20}-2 q^{22}-6 q^{24}-2 q^{26}\\
                   & -2 q^{28} + \cdots\\
\hat{Z}_{11} & = q^{\frac{1}{3}} \lb -q-q^3-3 q^5-2 q^7-3 q^9-2 q^{11}-5 q^{13}-2 q^{15}-4 q^{17}-2 q^{19}-5 q^{21}-4 q^{23}-2q^{25} \notag\right.\\
& \left. -q^{27} -3 q^{29}  + \cdots  \rb\\
\hat{Z}_{22} & = q^{\frac{1}{3}} \lb -q-q^3-3 q^5-2 q^7-3 q^9-2 q^{11}-5 q^{13}-2 q^{15}-4 q^{17}-2 q^{19}-5 q^{21}-4 q^{23}-2q^{25} \notag\right.\\
 & \left.  -q^{27} -3 q^{29}  + \cdots  \rb\\
\hat{Z}_{01} & = 3\zeta(-1,1/2) + \zeta(0,1/2) -2 q^4-q^6-2 q^8-3 q^{10}-3 q^{12}-2 q^{14}-4 q^{16}-2 q^{18}-5 q^{20} -3 q^{22}\\
&  -2 q^{24}-2 q^{26}-3q^{28} + \cdots \\
\hat{Z}_{02} & = 3\zeta(-1,1/2) + \zeta(0,1/2) -q^2-q^4-q^6-3 q^8-2 q^{10}-2 q^{12}-3 q^{14}-4 q^{16}-2 q^{18} -4 q^{20} \\
&  -2 q^{22} -3 q^{24}-2 q^{26}-4q^{28}  + \cdots \\
\hat{Z}_{12} & = q^{-\frac{1}{3}} \lb -2 q^3-2 q^5-2 q^7-2 q^9-4 q^{11}-2 q^{13}-4 q^{15}-2 q^{17}-6 q^{19}-2 q^{21}-4 q^{23}-2 q^{25} \notag\right.\\
& \left.   -4q^{27}   + \cdots  \rb  \\         
\end{align*}

\end{itemize}

\appendix
\section*{Appendix}
\addcontentsline{toc}{section}{Appendix}

\section{Type $\RN{1}$ Lie superalgebra and its quantization}

We give a summary of the representation theory of $osp(2|2n)$ and the quantum group $U_h(\text{type}\, \RN{1})$ in \cite{GP2}\footnote{For reviews on Lie superalgebras, see \cite{K1,K2,K3}}. For $osp(2|2n)$, the set of positive roots is $\Delta^{+}= \Delta^{+}_{\bar{0}} \cup \Delta^{+}_{\bar{1}}$ with
$$
\Delta^{+}_{\bar{0}} = \lac \delta_{i} \pm \delta_{j} | 1 \leq i < j \leq n \rac \cup \lac 2\delta_{i} \rac\quad \text{and}\quad \Delta^{+}_{\bar{1}}= \lac \epsilon \pm \delta_{i} \rac,\qquad n \in \intg_{+},
$$
where $\lac \epsilon, \delta_1,\cdots , \delta_n \rac$ form the dual basis of the Cartan subalgebra. Their inner products are
$$
\left\langle \epsilon, \epsilon \right\rangle =1 \qquad \left\langle \delta_i, \delta_j \right\rangle = - \delta_{ij}\qquad \left\langle \epsilon, \delta_j \right\rangle = 0
$$
The half sums of positive roots are
$$
\rho_{0} = \sum_{i} (n+1-i) \delta_{i},\quad \rho_{1}= n\epsilon \quad \text{and}\quad \rho= \rho_{0} - \rho_{1}.
$$
The fundamental weights are
$$
w_{1} = \epsilon\qquad w_{k+1} =  \epsilon + \sum_{i=1}^{k} \delta_{i}\quad \text{and}\quad k=1,\cdots,n.
$$
The weight $\lambda^{c}_{a}$ decomposes as
$$
\lambda^{c}_{a} = a w_{1} + c_{1}w_{2} + \cdots +  c_{n}w_{n+1},
$$
where $c \in \mathbb{N}^{r-1}, a\in \complex$
\newline

Let $g$ be a Lie superalgebra of type $\RN{1}$, $sl(m|n)$ or $osp(2|2n)\, (m\neq n)$.  $U_h(g)$ is the $\complex[[h]]$-Hopf superalgebra generated by
  $$
	E_{i}, F_{i}, h_{i},\quad i=1,\cdots , r=m+n-1\quad \text{or}\quad n+1
	$$
	satisfying the relations
	$$
	[ h_{i}, h_{j}] =0,\quad [ h_{i},E_{j}] = A_{ij}E_{j},\quad [ h_{i},F_{j}] = A_{ij}F_{j},
	$$
	$$
	[ E_{i},F_{j}] = \delta_{i,j}\frac{q^{h_i}- q^{-h_i}}{q-q^{-1}},\quad E_{s}^{2}= F_{s}^{2}=0
	$$
	and the quantum Serre-type relations; they are quadratic, cubic and quartic relations among $E_{i}$ or $F_{i}$ (see \cite{Yamane} Definition 4.2.1).  $A_{ij}$ is $r\times r$ Cartan matrix and $\lac s \rac =\tau \subset \lac 1,\cdots,r \rac$ determining the parity of the generators. $E_{s}, F_{s}$ are the only odd generators.  Moreover, the (anti)commutator is $[x,y]:= xy - (-1)^{\bar{x}\bar{y}}yx$.  The Hopf algebra structure is
	$$
	\Delta(E_{i})= \Delta(E_{i}) \otimes 1 + q^{-h_i} \otimes \Delta(E_{i}),\quad \epsilon(\Delta(E_{i}))=0,\quad S(\Delta(E_{i}))=-q^{h_i}\Delta(E_{i})
	$$
	$$
	\Delta(F_{i})= \Delta(F_{i}) \otimes 1 + q^{-h_i} \otimes \Delta(F_{i}),\quad \epsilon(\Delta(F_{i}))=0,\quad S(\Delta(E_{i}))=-\Delta(F_{i}) q^{h_i}
	$$
	$$
	\Delta(h_{i})= \Delta(F_{i}) \otimes 1 + 1 \otimes \Delta(h_{i}),\quad \epsilon(\Delta(h_{i}))=0,\quad S(\Delta(h_{i}))= -h_i
	$$

\section{Invariance checks}

We display the weighted positive definite generic plumbing graphs and their equivalent graphs of the examples in Section 4. Each vertex corresponds to a $S^1$-bundle over $S^2$; a positive integer is the Euler number of a bundle. The graphs are related by the Kirby-Neumman moves:
\begin{figure}[h!]
\centering
\includegraphics[scale=0.6]{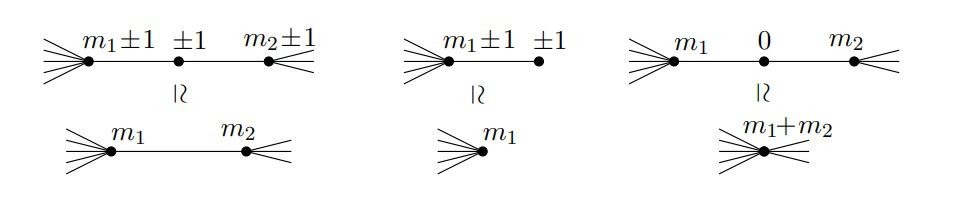}
\end{figure}
\newline
\noindent The invariance of $\hat{Z}_{ab}$ under the Kirby-Neumman moves for the exhibited graphs and the graphs mentioned in the captions of Figure 3 and 4 were verified. 
\begin{figure}[h!]
\begin{center}
\begin{tikzpicture}
\centering
\tikzstyle{every node}=[draw,shape=circle]

\draw (0,0) node[circle,fill,inner sep=1pt,label=left:$4$](4){} -- (1,0) node[circle,fill,inner sep=1pt,label=right:$1$](1){};
\draw (0,0)  -- (1,1) node[circle,fill,inner sep=1pt,label=right:$1$](1){};
\draw (0,0) -- (1,-1) node[circle,fill,inner sep=1pt,label=right:$1$](1){};

\end{tikzpicture}
\quad
\begin{tikzpicture}
\tikzstyle{every node}=[draw,shape=circle]

\draw (0,0) node[circle,fill,inner sep=1pt,label=left:$5$](5){} -- (1,0) node[circle,fill,inner sep=1pt,label=right:$1$](1){};
\draw (0,0)  -- (1,1) node[circle,fill,inner sep=1pt,label=right:$1$](1){};
\draw (0,0) -- (1,-1) node[circle,fill,inner sep=1pt,label=right:$1$](1){};
\draw (0,0) -- (0,1) node[circle,fill,inner sep=1pt,label=left:$1$](1){};

\end{tikzpicture}
\quad
\begin{tikzpicture}
\tikzstyle{every node}=[draw,shape=circle]

\draw (0,0)  -- (1,1) node[circle,fill,inner sep=1pt,label=above:$2$](2){} -- (2,1) node[circle,fill,inner sep=1pt,label=right:$1$](1){};
\draw (0,0) node[circle,fill,inner sep=1pt,label=left:$4$](4){} -- (1,0) node[circle,fill,inner sep=1pt,label=right:$1$](1){}  ;
\draw (0,0) -- (1,-1) node[circle,fill,inner sep=1pt,label=right:$1$](1){};

\end{tikzpicture}
\quad
\begin{tikzpicture}
\tikzstyle{every node}=[draw,shape=circle]

\draw (0,0) node[circle,fill,inner sep=1pt,label=left:$5$](5){} -- (1,0) node[circle,fill,inner sep=1pt,label=right:$1$](1){};
\draw (0,0)  -- (1,1) node[circle,fill,inner sep=1pt,label=right:$1$](1){};
\draw (0,0) -- (1,-1) node[circle,fill,inner sep=1pt,label=below:$1$](1){}-- (2,-1) node[circle,fill,inner sep=1pt,label=right:$1$](1){};

\end{tikzpicture}
\quad
\begin{tikzpicture}
\tikzstyle{every node}=[draw,shape=circle]

\draw (0,0) node[circle,fill,inner sep=1pt,label=above:$1$](1){} -- (1,0) node[circle,fill,inner sep=1pt,label=above:$2$](2){} -- (2,0) node[circle,fill,inner sep=1pt,label=above:$0$](0){} -- (3,0) node[circle,fill,inner sep=1pt,label=above:$2$](2){} -- (4,1) node[circle,fill,inner sep=1pt,label=right:$1$](1){};
\draw (3,0) node[circle,fill,inner sep=1pt,label=above:$2$](2){} -- (4,-1) node[circle,fill,inner sep=1pt,label=right:$1$](1){};

\end{tikzpicture}

\end{center}
\caption{A generic plumbing graph of $S^3$ (left most) and its equivalent graphs}
\end{figure}
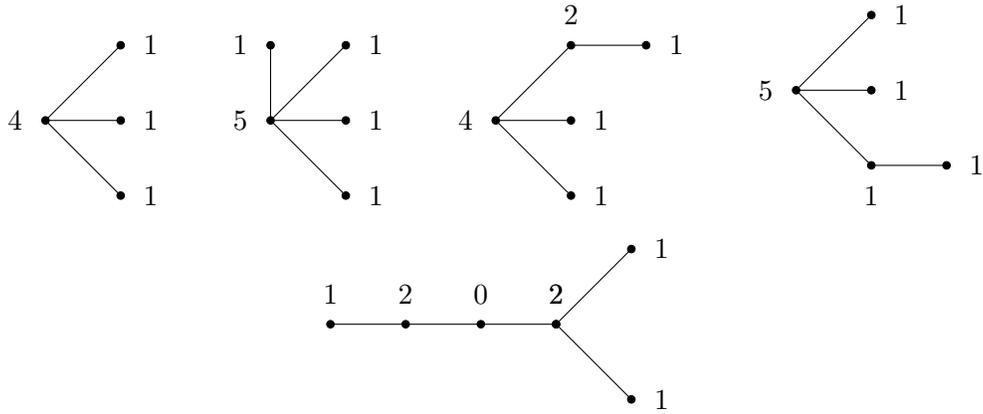

\begin{figure}[h!]
\begin{center}
\begin{tikzpicture}
\centering
\tikzstyle{every node}=[draw,shape=circle]

\draw (0,0) node[circle,fill,inner sep=1pt,label=left:$1$](1){} -- (1,0) node[circle,fill,inner sep=1pt,label=right:$2$](2){};
\draw (0,0)  -- (1,1) node[circle,fill,inner sep=1pt,label=right:$3$](3){};
\draw (0,0) -- (1,-1) node[circle,fill,inner sep=1pt,label=right:$7$](7){};

\end{tikzpicture}
\quad
\begin{tikzpicture}
\tikzstyle{every node}=[draw,shape=circle]

\draw (0,0) node[circle,fill,inner sep=1pt,label=left:$2$](2){} -- (1,0) node[circle,fill,inner sep=1pt,label=right:$3$](3){};
\draw (0,0) -- (0,1) node[circle,fill,inner sep=1pt,label=left:$1$](1){};
\draw (0,0)  -- (1,1) node[circle,fill,inner sep=1pt,label=right:$2$](2){};
\draw (0,0) -- (1,-1) node[circle,fill,inner sep=1pt,label=right:$7$](7){};

\end{tikzpicture}
\quad
\begin{tikzpicture}
\tikzstyle{every node}=[draw,shape=circle]

\draw (0,0) node[circle,fill,inner sep=1pt,label=left:$1$](1){} -- (1,0) node[circle,fill,inner sep=1pt,label=right:$3$](3){}  ;
\draw (0,0)  -- (1,1) node[circle,fill,inner sep=1pt,label=above:$3$](3){} -- (2,1) node[circle,fill,inner sep=1pt,label=right:$1$](1){};
\draw (0,0) -- (1,-1) node[circle,fill,inner sep=1pt,label=right:$7$](7){};

\end{tikzpicture}
\quad
\begin{tikzpicture}
\tikzstyle{every node}=[draw,shape=circle]

\draw (0,0) node[circle,fill,inner sep=1pt,label=left:$1$](1){} -- (1,0) node[circle,fill,inner sep=1pt,label=above:$4$](4){}-- (2,0) node[circle,fill,inner sep=1pt,label=right:$1$](1){}   ;
\draw (0,0)  -- (1,1) node[circle,fill,inner sep=1pt,label=right:$2$](2){};
\draw (0,0) -- (1,-1) node[circle,fill,inner sep=1pt,label=right:$7$](7){};

\end{tikzpicture}
\quad
\begin{tikzpicture}
\tikzstyle{every node}=[draw,shape=circle]

\draw (0,0) node[circle,fill,inner sep=1pt,label=left:$1$](1){} -- (1,0) node[circle,fill,inner sep=1pt,label=right:$3$](3){};
\draw (0,0)  -- (1,1) node[circle,fill,inner sep=1pt,label=right:$2$](2){};
\draw (0,0) -- (1,-1) node[circle,fill,inner sep=1pt,label=above:$8$](8){}-- (2,-1) node[circle,fill,inner sep=1pt,label=right:$1$](1){};

\end{tikzpicture}
\quad
\begin{tikzpicture}
\tikzstyle{every node}=[draw,shape=circle]

\draw (0,0) node[circle,fill,inner sep=1pt,label=left:$2$](2){} -- (1,0) node[circle,fill,inner sep=1pt,label=right:$3$](3){}  ;
\draw (0,0)  -- (1,1) node[circle,fill,inner sep=1pt,label=above:$1$](1){} -- (2,1) node[circle,fill,inner sep=1pt,label=right:$3$](3){};
\draw (0,0) -- (1,-1) node[circle,fill,inner sep=1pt,label=right:$7$](7){};

\end{tikzpicture}
\quad
\begin{tikzpicture}
\tikzstyle{every node}=[draw,shape=circle]

\draw (0,0) node[circle,fill,inner sep=1pt,label=left:$2$](2){} -- (1,0) node[circle,fill,inner sep=1pt,label=above:$1$](1){}-- (2,0) node[circle,fill,inner sep=1pt,label=right:$4$](4){}   ;
\draw (0,0)  -- (1,1) node[circle,fill,inner sep=1pt,label=right:$2$](2){};
\draw (0,0) -- (1,-1) node[circle,fill,inner sep=1pt,label=right:$7$](7){};

\end{tikzpicture}
\quad
\begin{tikzpicture}
\tikzstyle{every node}=[draw,shape=circle]

\draw (0,0) node[circle,fill,inner sep=1pt,label=left:$2$](2){} -- (1,0) node[circle,fill,inner sep=1pt,label=right:$3$](3){};
\draw (0,0)  -- (1,1) node[circle,fill,inner sep=1pt,label=right:$2$](2){};
\draw (0,0) -- (1,-1) node[circle,fill,inner sep=1pt,label=above:$1$](1){}-- (2,-1) node[circle,fill,inner sep=1pt,label=right:$8$](8){};

\end{tikzpicture}
\quad
\begin{tikzpicture}
\tikzstyle{every node}=[draw,shape=circle]

\draw (0,0) node[circle,fill,inner sep=1pt,label=above:$2$](2){} -- (1,0) node[circle,fill,inner sep=1pt,label=above:$1$](1){} -- (2,0) node[circle,fill,inner sep=1pt,label=above:$0$](0){} -- (3,0) node[circle,fill,inner sep=1pt,label=above:$0$](0){} -- (4,1) node[circle,fill,inner sep=1pt,label=right:$3$](3){};
\draw (3,0) node[circle,fill,inner sep=1pt,label=above:$0$](0){} -- (4,-1) node[circle,fill,inner sep=1pt,label=right:$7$](7){};

\end{tikzpicture}
\quad
\begin{tikzpicture}
\tikzstyle{every node}=[draw,shape=circle]

\draw (0,0) node[circle,fill,inner sep=1pt,label=above:$3$](3){} -- (1,0) node[circle,fill,inner sep=1pt,label=above:$1$](1){} -- (2,0) node[circle,fill,inner sep=1pt,label=above:$0$](0){} -- (3,0) node[circle,fill,inner sep=1pt,label=above:$0$](0){} -- (4,1) node[circle,fill,inner sep=1pt,label=right:$2$](2){};
\draw (3,0) node[circle,fill,inner sep=1pt,label=above:$0$](0){} -- (4,-1) node[circle,fill,inner sep=1pt,label=right:$7$](7){};

\end{tikzpicture}
\quad
\begin{tikzpicture}
\tikzstyle{every node}=[draw,shape=circle]

\draw (0,0) node[circle,fill,inner sep=1pt,label=above:$7$](7){} -- (1,0) node[circle,fill,inner sep=1pt,label=above:$1$](1){} -- (2,0) node[circle,fill,inner sep=1pt,label=above:$0$](0){} -- (3,0) node[circle,fill,inner sep=1pt,label=above:$0$](0){} -- (4,1) node[circle,fill,inner sep=1pt,label=right:$2$](2){};
\draw (3,0) node[circle,fill,inner sep=1pt,label=above:$0$](0){} -- (4,-1) node[circle,fill,inner sep=1pt,label=right:$3$](3){};

\end{tikzpicture}
\end{center}
\caption{A generic plumbing graph of $\Sigma(2,3,7)$ (the first graph) and its equivalent graphs} 
\end{figure}
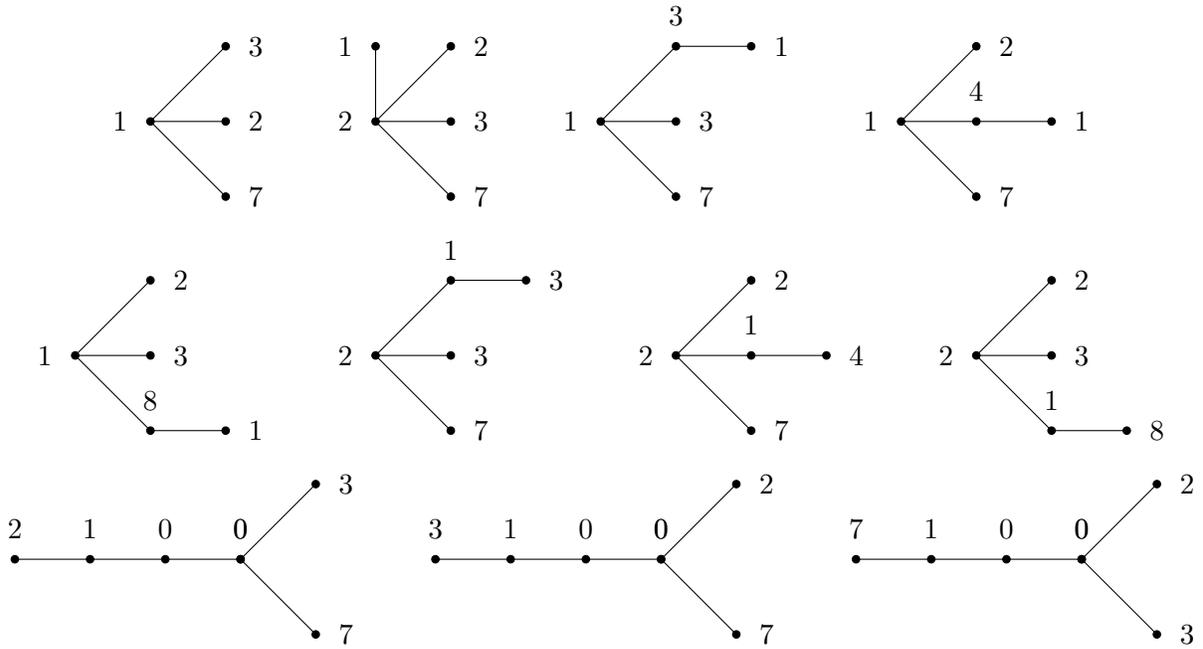

\clearpage

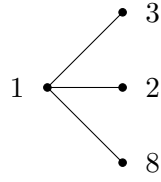
\begin{figure}[h]
\begin{center}
\begin{tikzpicture}
\centering
\tikzstyle{every node}=[draw,shape=circle]

\draw (0,0) node[circle,fill,inner sep=1pt,label=left:$1$](1){} -- (1,0) node[circle,fill,inner sep=1pt,label=right:$2$](2){};
\draw (0,0)  -- (1,1) node[circle,fill,inner sep=1pt,label=right:$3$](3){};
\draw (0,0) -- (1,-1) node[circle,fill,inner sep=1pt,label=right:$8$](8){};

\end{tikzpicture}
\end{center}
\caption{A generic plumbing graph of $M(-1 \vert \frac{1}{2}, \frac{1}{3}, \frac{1}{8})$ and its equivalent graphs can be obtained from Figure 2 by replacing the weight 7 by 8 or 8 by 9.}
\end{figure}

\begin{figure}[h]
\begin{center}
\begin{tikzpicture}
\centering
\tikzstyle{every node}=[draw,shape=circle]

\draw (0,0) node[circle,fill,inner sep=1pt,label=left:$1$](1){} -- (1,0) node[circle,fill,inner sep=1pt,label=right:$2$](2){};
\draw (0,0)  -- (1,1) node[circle,fill,inner sep=1pt,label=right:$3$](3){};
\draw (0,0) -- (1,-1) node[circle,fill,inner sep=1pt,label=right:$9$](9){};

\end{tikzpicture}
\end{center}
\caption{A generic plumbing graph of $M(-1 \vert \frac{1}{2}, \frac{1}{3}, \frac{1}{9})$ and its equivalent graphs can be obtained from Figure 2 by replacing the weight 7 by 9 or 8 by 10.}
\end{figure}
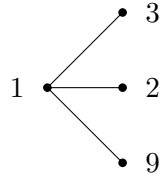

\section{Derivations}

We setup for the derivations of the ingredients in the CGP invariant formula using the appendix A and outline the derivation of $\hat{Z}^{osp(2|2)}$ in Section 3.
\newline

The root lattice $\Lambda_R$ of $osp(2|2)$ are generated by $2\delta$ and $\epsilon - \delta$, hence, two dimensional. The pivotal element $\pi$
\begin{equation*}
\begin{split}
\pi & = 2 ( \rho_{0} - \rho_{1} ) - 2l\rho_{0} \in \Lambda_R\\
    & = -2 (\epsilon - \delta) - l (2\delta)
\end{split}
\end{equation*}
This implies that 
$$
K_{\pi} = K_{1}^{-l} K_{2}^{-2} \in U^{H}_{\xi}(osp(2|2))
$$
The weight $\lambda$ of $V$ is
\begin{equation*}
\begin{split}
\lambda & = \mu_1 w_1 + \mu_2 w_2 \\
       & =  ( \mu_1 + \mu_2 ) \epsilon + \mu_2 \delta,
\end{split}
\end{equation*}
where $w_1= \epsilon$ and $w_2=  \epsilon + \delta$. Under the assumption mentioned in Section 2, we arrive at 
$$
\theta_V (\vec{\mu}; l) = \xi^{2l\mu_{2}}\xi^{\mu_{1}^2 + 2\mu_{1}\mu_{2}-4\mu_{2}-2\mu_{1}}\, 1_V\qquad l= \text{odd and} \geq 3.  
$$
\begin{equation*}
\begin{split}
S(\vec{\mu},\vec{\mu}^{\prime}; l) & = \xi^{2l(\mu_{2}+\mu^{\prime}_{2})} \xi^{2 \lb \mu_{1}\mu^{\prime}_{1}+\mu_{1}\mu^{\prime}_{2}+\mu^{\prime}_{1}\mu_{2} \rb-2 \lb \mu_{1} + \mu^{\prime}_{1} + 2 \mu_{2}+ 2\mu^{\prime}_{2} \rb}\\
& \times \frac{\lac 2l(\mu^{\prime}_{2} +1 -l) \rac}{\lac 2(\mu^{\prime}_{2} +1 -l) \rac} \lac \mu^{\prime}_{1} -2 +l \rac \lac  \mu^{\prime}_{1} + 2\mu^{\prime}_{2} - l \rac
\end{split}
\end{equation*}
$$
d(\vec{\mu}; l) = \frac{\lac 2(\mu_{2} +1 -l) \rac}{\lac 2l(\mu_{2} +1 -l) \rac} \frac{1}{\lac \mu_{1} -2 +l \rac \lac  \mu_{1} + 2\mu_{2} - l \rac},
$$
$$
\lac x \rac_{\xi} := \xi^x - \xi^{-x}\quad \xi=q^{1/2}=e^{i 2\pi/l},
$$
where $(\mu_{1}, \mu_{2})$ and $(\mu^{\prime}_{1}, \mu^{\prime}_{2})$ are the coefficients in the weight decompositions of $\lambda$ and $\mu$ for $V$ and $V^{\prime}$, respectively (see appendix A). For the S-matrix, notations for $V$ and $V^{\prime}$ are switched.  
\newline

In order to apply the derivation of the superalgebra $\hat{Z}$ given in \cite{FP}, we need to adapt the above three expressions into a plumbing graph setup. We first let
$$
\vec{\alpha} = ( \alpha_1 = \mu_1 -2 + l,\, \alpha_2 = 2\mu_2 + 2 -2l)\quad \in \complex^2
$$
Then 
$$
\theta_{V} = \xi^{\alpha_{1}^{2}+\alpha_{1}\alpha_{2}}
$$
$$
S^{\prime} = \xi^{2\alpha_{1}\alpha^{\prime}_{1} + \alpha_{1}\alpha^{\prime}_{2} +\alpha_{2}\alpha^{\prime}_{1}} \frac{\lac l \alpha^{\prime}_{2}\rac_{\xi} \lac \alpha^{\prime}_{1}\rac_{\xi} \lac \alpha^{\prime}_{1} + \alpha^{\prime}_{2}\rac_{\xi}}{\lac \alpha^{\prime}_{2} \rac_{\xi}}
$$
We next shift $\alpha_1$ and $\alpha_2$ by s and t, respectively. And then we associate $\theta_V$ to each vertex. This leads to
$$
\theta_{V_{\alpha^{I}_{s^{I} t^{I}}}} = \xi^{(\mu_1 -2 +s)^{I}(\mu_1 + 2\mu_2 +s + t)^{I}}.
$$
Similarly, the S-matrix elements between to vertices I and J of graphs contain
$$
\prod_{(I,J)\in Edges} S^{\prime}(\alpha^{J}_{s^{J} t^{J}}, \alpha^{I}_{s^{I} t^{I}} ) \supset \xi^{\sum_{IJ }B_{IJ}(\mu_1 -2 +s)^{I}(\mu_1 + 2\mu_2 +s + t)^{J} },
$$
which is the relevant part for the derivation. For the modified quantum dimension $d(\vec{\alpha})$
$$
d(\vec{\alpha}) = \frac{\lac \alpha_{2} \rac_{\xi}}{\lac l \alpha_{2}\rac_{\xi} \lac \alpha_{1}\rac_{\xi} \lac \alpha_{1} + \alpha_{2}\rac_{\xi}},
$$
after shifting by s and t as above and some manipulations, the roots of unity dependent factor becomes
$$
d(\vec{\alpha} ) \supset \frac{ \xi^{2(\mu_{1} + 2\mu_{2}  + s+t )} - \xi^{2(\mu_{1} + s-2 )} }{  \lb 1- \xi^{2(\mu_{1} + 2\mu_{2}  + s+t )} \rb \lb 1- \xi^{2(\mu_{1} + s-2 )} \rb} 
$$
We define coordinates of the 2D Cartan subalgebra of $osp(2|2)$ to be  
$$
y= \xi^{2(\mu_{1} + 2\mu_{2}  + s+t )} \qquad z= \xi^{2(\mu_{1} + s-2 )}
$$
Then the modified quantum dimension for a graph contains 
$$
d(y,z) \supset \frac{y_{I} -z_{I}}{(1-y_{I})(1-z_{I})}
$$
\newline

We next sketch the derivation of $\hat{Z}^{osp(2|2)}$ in Section 3 by applying the procedure in the appendix D of \cite{FP}.  For a closed oriented graph 3-manifold $Y$ equipped with $\omega \in H^{1}(Y; \rational / \intg \times  \rational / \intg )$, the CGP invariant for a pair $(Y,\omega)$ in which Y is presented by Dehn surgery is  
$$
N_{l}(Y,\omega)  = \sum_{s^{I}, t^{I}=0 }^{l-1} d(\alpha^{I_{0}}_{s^{I_{0}} t^{I_{0}}}) \prod_{I\in Vert} d(\alpha^{I}_{s^{I} t^{I}}) \left\langle \theta_{V_{\alpha^{I}_{s^{I} t^{I}}}} \right\rangle^{B_{II}} \prod_{(I,J)\in Edges} S^{\prime}(\alpha^{J}_{s^{J} t^{J}}, \alpha^{I}_{s^{I} t^{I}} )
$$
$$
\omega \in H^{1}(Y; \rational / \intg \times  \rational / \intg ) \cong B^{-1}\intg^L / \intg^L \times B^{-1}\intg^L / \intg^L
$$
$$
\omega([m_I]) = \mu^{I} = ( \mu^{I}_{1} , \mu^{I}_{2} ) \in \rational / \intg \times  \rational / \intg,\qquad \sum_{J} B_{IJ}\mu^{J} = 0\quad \text{mod}\, \intg \times \intg,
$$
where $m_I$ is a meridian of I-th component $L_I$ of a surgery link $L$ and $[m_I]$ is its homology class and $B$ is a linking matrix of $L$. A color of $L_I$ is set by $\omega([m_I])$. The origin of $\omega$ can be found in the footnote in Section 1.3 of \cite{CGP}. After substituting the ingredients, the right hand side becomes
\begin{equation*}
\begin{split}
N_{l}(Y,\omega)   & = \frac{1}{l^{L+1}} \prod_{I\in V} \lb e^{i2\pi \mu^{I}_{1}} - e^{-i2\pi \mu^{I}_{1} } \rb^{deg(I)-2}\\
	& \times \sum_{a^{I}, b^{I} \in \intg^{L}/ l\intg^{L} } F \lb \lac \xi^{2(\mu_{1} + 2\mu_{2}  + a )} ,  \xi^{2(\mu_{1} + b )} \rac_{I\in V} \rb \xi^{\sum_{IJ}B_{IJ} (\mu_{1} + 2\mu_{2}  + a )^{I} (\mu_{1} + b )^{J} },
\end{split}
\end{equation*}
where $a^{I}=s^{I} +t^{I},\, b^{I}=s^{I}-2$ and 
$$
F(y,z) = \prod_{I\in V} \lb \frac{y_{I} - z_{I}}{(1-y_{I})(1-z_{I})} \rb^{2- deg(I)}.
$$
We next expand $F(y,z)$, which modifies the summand as
$$
\sum_{a^{I}, b^{I} \in \intg^{L}/ l\intg^{L} }  \xi^{\sum_{IJ} B_{IJ} (\mu_{1} + 2\mu_{2}  + a )^{I} (\mu_{1} + b )^{J} +  2\sum_{I} n_{I}(\mu_{1} + 2\mu_{2}  + a)^{I} + m_{I}(\mu_{1} + b )^{I} }
$$
And then we recast it in terms of matrices by forming
$$
\mathcal{M} = \frac{1}{2}\begin{pmatrix}
0 & B \\
B & 0 
\end{pmatrix}\qquad
r= \begin{pmatrix}
b \\
a
\end{pmatrix}\qquad
p= \begin{pmatrix}
B (\mu_{1} + 2\mu_{2}) + 2m \\
B \mu_{1} + 2n
\end{pmatrix},
$$
This enables us to apply the appropriate Gauss reciprocity formula 
\begin{multline*}
\sum_{r \in \intg^L/ l\intg^L} exp \lb \frac{i2\pi}{l} r^T \mathcal{M} r + \frac{i2\pi}{l} p^T r \rb = \\
\frac{e^{i\pi \sigma(\mathcal{M})/4} ( l/2)^{N/2}}{| Det \mathcal{M}|^{1/2} } \sum_{\delta \in \intg^L/ 2\mathcal{M}\intg^L} exp \lb -\frac{i\pi l}{2} \lb \delta + \frac{p}{l} \rb^T \mathcal{M}^{-1} \lb \delta + \frac{p}{l} \rb \rb,
\end{multline*} 
where $\mathcal{M}$ is a non-degenerate symmetric  $2L \times 2L$ matrix over $\intg$ and $\sigma(\mathcal{M})$ is its signature. From the p-quadratic term in the exponential we read off 
$$
RHS \supset \xi^{-4 m^{T}B^{-1}n},
$$ 
which becomes the q-term in the summand in Section 3 after we analytically continue from a complex unit circle into an interior of an unit disc coordinatized by q.

\end{document}